\documentclass[11pt]{article}

\makeatletter

\@addtoreset{equation}{section}
\makeatother

\newtheorem{definition}{Definition}[section]
\newtheorem{theorem}{Theorem}[section]

\newtheorem{remark}{Remark}[section]

\usepackage{amssymb}
\usepackage{latexsym}

\begin{document} 

\title{Bounds on Non-Hermite Wigner-Yanase-Dyson Skew Information} 
\author{Kenjiro Yanagi \\
Emeritus Professor, Yamaguchi University \\
{\tt yanagi@yamaguchi-u.ac.jp}}
\date{}

\maketitle


\begin{abstract}
We obtain the relation between non-hermite skew information $NI_{\rho,\alpha}(A)$, $NU_{\rho,\alpha}(A)$, variance $NV_{\rho}(A)$ and $NJ_{\rho,\alpha}(A)$ 
by using the maximal eigenvalue and the minimal eigenvalue of density $\rho$. And also we obtain their refinements. 

\end{abstract}


\section{Introduction}

In quantum mechanical system, the expectation value of a non-hermite  observable $A$ in a quantum state $\rho$ 
is expressed by $E_{\rho}(A) = Tr[\rho A]$.  It is natural that the variance for a quantum state $\rho$ and an observable 
$A$ is defined by 
$$
NV_{\rho}(A) = \frac{1}{2}Tr[\rho(A_0^*A_0+A_0A_0^*)],
$$
where $A_0 = A - (Tr[\rho A])I$. 
It is well known that non-hermite Robertson uncertainty relation is represented by the commutator $[A,B] = AB-BA$ 
in the following inequality.
$$
\frac{1}{4} |Tr[\rho [A,B]]|^2 \leq NV_{\rho}(A) NV_{\rho}(B).
$$

\noindent
The different non-hermite Wigner-Yanase-Dyson type of uncertainty relation is given as follows:

\begin{definition}
Let $\rho$ be a density operator, $A, B$ be non-hermite matrices and $0 \leq \alpha \leq 1$. 
We define the following non-hermite Wigner-Yanas-Dyson skew information $NI_{\rho,\alpha}(A)$,  $NJ_{\rho,\alpha}(A)$ and $NU_{\rho,\alpha}(A)$.
\begin{eqnarray*}
NI_{\rho,\alpha}(A) & = & \frac{1}{2} Tr[(i[\rho^{\alpha},A_0^*])(i[\rho^{1-\alpha},A_0])] \\
& = & \frac{1}{2}Tr[\rho (A_0^*A_0+A_0A_0^*)] - \frac{1}{2}Tr[\rho^{\alpha}A_0^* \rho^{1-\alpha}A_0+ \rho^{\alpha}A_0\rho^{1-\alpha}A_0^*]
\end{eqnarray*}
\begin{eqnarray*}
NJ_{\rho,\alpha}(A) & = & \frac{1}{2} Tr[ \{ \rho^{\alpha}, A_0^* \} \{ \rho^{1-\alpha},A_0 \}] \\
& = & \frac{1}{2}Tr[\rho (A_0^*A_0+A_0A_0^*)] + \frac{1}{2}Tr[\rho^{\alpha}A_0^* \rho^{1-\alpha}A_0+ \rho^{\alpha}A_0\rho^{1-\alpha}A_0^*]
\end{eqnarray*}
$$
NU_{\rho,\alpha}(A) = \sqrt{ NI_{\rho,\alpha}(A) NJ_{\rho,\alpha}(A)},
$$
where $\{ A,B \} = AB+BA$ is anti-commutator.  \\
We remark that $0 \leq NI_{\rho,\alpha}(A) \leq NU_{\rho,\alpha}(A) \leq NV_{\rho}(A) \leq NJ_{\rho,\alpha}(A)$.
\label{def:definition1.1}
\end{definition}

\begin{theorem} [\cite{Ya:ur-wynh}]
Let $\rho$ be a density operator and $A, B$ be non-hermite matrices. .  
Then the following uncertainty relation holds.
$$
\alpha(1-\alpha) |Tr[\rho [A,B]]|^2 \leq NU_{\rho,\alpha}(A) NU_{\rho,\alpha}(B),  
$$
\label{th:theorem1.1}
\end{theorem}

\noindent
Proof of Theorem \ref{th:theorem1.1}~~Let $\rho = \sum_{i,j} \lambda_i |\phi_i \rangle \langle \phi_i|$ 
be the spectral decomposition. Since 
$$
NI_{\rho,\alpha}(A) = \frac{1}{2} \sum_{i,j} (\lambda_i^{\alpha}-\lambda_j^{\alpha})(\lambda_i^{1-\alpha}-\lambda_j^{1-\alpha}) |a_{ij}|^2,
$$
$$
NJ_{\rho,\alpha}(A) = \frac{1}{2} \sum_{i,j} (\lambda_i^{\alpha}+\lambda_j^{\alpha})(\lambda_i^{1-\alpha}+\lambda_j^{1-\alpha}) |a_{ij}|^2, 
$$
where $a_{ij} = \langle \phi_i|A_0|\phi_j \rangle$, we have 
\begin{eqnarray*}
&   & \alpha(1-\alpha) |Tr[\rho[A,B]]|^2 \\
& = & \alpha(1-\alpha) |\sum_{i,j} (\lambda_i-\lambda_j) a_{ij} b_{ji}|^2 \\
& \leq & \frac{1}{4} \{  \sum_{i,j} 2\sqrt{\alpha(1-\alpha)} |\lambda_i-\lambda_j| |a_{ij}||b_{ji}| \}^2 \\
& \leq & \frac{1}{4} \{ \sum_{i,j} \{ ((\lambda_i+\lambda_j)^2-(\lambda_i^{\alpha}\lambda_j^{1-\alpha}+\lambda_i^{1-\alpha}\lambda_j^{\alpha})^2 \}^{1/2} |a_{ij}| |b_{ji}| \}^2 \\
& \leq & \frac{1}{2} \sum_{i,j} (\lambda_i+\lambda_j-\lambda_i^{\alpha}\lambda_j^{1-\alpha}-\lambda_i^{1-\alpha}\lambda_j^{\alpha})|a_{ij}|^2~\frac{1}{2} \sum_{i,j} (\lambda_i+\lambda_j+\lambda_i^{\alpha}\lambda_j^{1-\alpha}+\lambda_i^{1+\alpha}\lambda_j^{\alpha})|b_{ji}|^2 \\
& = & NI_{\rho,\alpha}(A) NJ_{\rho,\alpha}(B).
\end{eqnarray*}
Since $4\alpha(1-\alpha)(t-1)^2 \leq (t+1)^2-(t^{\alpha}+t^{1-\alpha})^2$ for any $t > 0$ and $0 \leq \alpha \leq 1$, 
the second inequality holds. We also have 
$$
\alpha(1-\alpha)|Tr[\rho[A,B]]|^2 \leq NI_{\rho,\alpha}(B) NJ_{\rho,\alpha}(A).
$$
Hence we have the final result. \ \hfill $\Box$

\vspace{0.5cm}
\noindent
In this paper we obtain the refined relation between $NI_{\rho,\alpha}(A)$, $NU_{\rho,\alpha}(A)$, $NV_{\rho}(A)$ and $NJ_{\rho,\alpha}(A)$ by using 
the maximal eigenvalue and minimal eigenvalue of $\rho$ in section 2. 
And also we obtain their refinements in section 3. 


\section{Non-Hermite Wigner-Yanase-Dyson Skew Information}

\begin{theorem}
Let $\rho$ be a density matrix with eigenvalues $0 < \lambda_1 \leq \lambda_2 \leq \cdots \leq \lambda_n$ and $\lambda_1 \neq \lambda_n$ and let $0 \leq \alpha \leq 1$. Then
\begin{description}
\item[(1)]  $\displaystyle{ \frac{\lambda_n+\lambda_1}{(\lambda_n^{\alpha}-\lambda_1^{\alpha})(\lambda_n^{1-\alpha}-\lambda_1^{1-\alpha})} NI_{\rho,\alpha}(A) \leq NV_{\rho}(A) \leq \frac{\lambda_n+\lambda_1}{(\lambda_n^{\alpha}+\lambda_1^{\alpha})(\lambda_n^{1-\alpha}+\lambda_1^{1-\alpha})} NJ_{\rho,\alpha}(A)}$.
\item[(2)]  $\displaystyle{ \left(\frac{(\lambda_n^{\alpha}+\lambda_1^{\alpha})(\lambda_n^{1-\alpha}+\lambda_1^{1-\alpha})}{(\lambda_n^{\alpha}-\lambda_1^{\alpha})(\lambda_n^{1-\alpha}-\lambda_1^{1-\alpha})} \right)^{1/2} NI_{\rho,\alpha}(A) }$ \\
$\displaystyle{ \leq NU_{\rho,\alpha}(A) \leq \left( \frac{(\lambda_n^{\alpha}-\lambda_1^{\alpha})(\lambda_n^{1-\alpha}-\lambda_1^{1-\alpha})}{(\lambda_n^{\alpha}+\lambda_1^{\alpha})(\lambda_n^{1-\alpha}+\lambda_1^{1-\alpha})} \right)^{1/2} NJ_{\rho,\alpha}(A) }$.
\item[(3)]  $\displaystyle{ NU_{\rho,\alpha}(A) \leq \frac{ \{ (\lambda_n^{\alpha}-\lambda_1^{\alpha})(\lambda_n^{\alpha}+\lambda_1^{\alpha})(\lambda_n^{1-\alpha}-\lambda_1^{1-\alpha})(\lambda_n^{1-\alpha}+\lambda_1^{1-\alpha}) \}^{1/2}}{\lambda_n+\lambda_1} NV_{\rho}(A) }$.
\end{description}
\label{th:theorem2.1}
\end{theorem}

\noindent
Proof of Theorem \ref{th:theorem2.1}. 
(1)~~We show the relation between $NI_{\rho,\alpha}(A)$, $NV_{\rho}(A)$ and $NJ_{\rho,\alpha}(A)$. 
Since 
$$
\max \left\{ \frac{(t^{\alpha}-1)(t^{1-\alpha}-1)}{t+1} ; \frac{\lambda_1}{\lambda_n} \leq t \leq \frac{\lambda_n}{\lambda_1} \right\} = \frac{(\lambda_n^{\alpha} - \lambda_1^{\alpha})(\lambda_n^{1-\alpha}-\lambda_1^{1-\alpha})}{\lambda_n+\lambda_1},
$$
we have
\begin{eqnarray}
NI_{\rho,\alpha}(A) & = & \frac{1}{2} \sum_{i,j} (\lambda_i^{\alpha}-\lambda_j^{\alpha})(\lambda_i^{1-\alpha}-\lambda_j^{1-\alpha}) |a_{ij}|^2 \nonumber \\
& = & \frac{1}{2} \sum_{i,j} \frac{ (\lambda_i^{\alpha}-\lambda_j^{\alpha})(\lambda_i^{1-\alpha}-\lambda_j^{1-\alpha})}{\lambda_i+\lambda_j} (\lambda_i+\lambda_j) |a_{ij}|^2 \nonumber \\
& \leq & \max_{i,j} \frac{ (\lambda_i^{\alpha}-\lambda_j^{\alpha})(\lambda_i^{1-\alpha}-\lambda_j^{1-\alpha})}{\lambda_i+\lambda_j} \frac{1}{2} \sum_{i,j} (\lambda_i+\lambda_j)|a_{ij}|^2 \nonumber \\
& = & \frac{(\lambda_n^{\alpha}-\lambda_1^{\alpha})(\lambda_n^{1-\alpha}-\lambda_1^{1-\alpha})}{\lambda_n+\lambda_1} NV_{\rho}(A).  \label{eq:num2.1}
\end{eqnarray}
Since
$$
\max \left\{ \frac{t+1}{(t^{\alpha}+1)(t^{1-\alpha}+1)} ; \frac{\lambda_1}{\lambda_n} \leq t \leq \frac{\lambda_n}{\lambda_1} \right\} = \frac{\lambda_n +\lambda_1}{(\lambda_n^{\alpha}+\lambda_1^{\alpha})(\lambda_n^{1-\alpha}+\lambda_1^{1-\alpha})}, 
$$
we have 
\begin{eqnarray}
NV_{\rho}(A) & = & \frac{1}{2} \sum_{i,j} (\lambda_i+\lambda_j) |a_{ij}|^2 \nonumber \\
& = & \frac{1}{2} \sum_{i,j} \frac{\lambda_i+\lambda_j}{(\lambda_i^{\alpha}+\lambda_j^{\alpha})(\lambda_i^{1-\alpha}+\lambda_j^{1-\alpha})}
(\lambda_i^{\alpha}+\lambda_j^{\alpha})(\lambda_i^{1-\alpha}+\lambda_j^{1-\alpha}) |a_{ij}|^2 \nonumber \\
& \leq & \max_{i,j} \frac{\lambda_i+\lambda_j}{(\lambda_i^{\alpha}+\lambda_j^{\alpha})(\lambda_i^{1-\alpha}+\lambda_j^{1-\alpha})} \frac{1}{2} \sum_{i,j} (\lambda_i^{\alpha}+\lambda_j^{\alpha})(\lambda_i^{1-\alpha}+\lambda_j^{1-\alpha}) |a_{ij}|^2 \nonumber \\
& = &  \frac{\lambda_n +\lambda_1}{(\lambda_n^{\alpha}+\lambda_1^{\alpha})(\lambda_n^{1-\alpha}+\lambda_1^{1-\alpha})}NJ_{\rho,\alpha}(A). \label{eq:num2.2}
\end{eqnarray}
Then we have 
$$
\frac{\lambda_n+\lambda_1}{(\lambda_n^{\alpha}-\lambda_1^{\alpha})(\lambda_n^{1-\alpha}-\lambda_1^{1-\alpha})} NI_{\rho,\alpha}(A) \leq NV_{\rho}(A) \leq  \frac{\lambda_n +\lambda_1}{(\lambda_n^{\alpha}+\lambda_1^{\alpha})(\lambda_n^{1-\alpha}+\lambda_1^{1-\alpha})} NJ_{\rho}(A).
$$
(2)~~We show the relation between $NI_{\rho,\alpha}(A)$, $NU_{\rho,\alpha}(A)$ and $NJ_{\rho,\alpha}(A)$. By (\ref{eq:num2.1}) and (\ref{eq:num2.2}) we have 
$$
NI_{\rho,\alpha}(A) \leq \frac{(\lambda_n^{\alpha}-\lambda_1^{\alpha})(\lambda_n^{1-\alpha}-\lambda_1^{1-\alpha})}{(\lambda_n^{\alpha}+\lambda_1^{\alpha})(\lambda_n^{1-\alpha}+\lambda_1^{1-\alpha})} NJ_{\rho,\alpha}(A).
$$
Since 
$$
NI_{\rho,\alpha}(A)^2 \leq \frac{(\lambda_n^{\alpha}-\lambda_1^{\alpha})(\lambda_n^{1-\alpha}-\lambda_1^{1-\alpha})}{(\lambda_n^{\alpha}+\lambda_1^{\alpha})(\lambda_n^{1-\alpha}+\lambda_1^{1-\alpha})} NI_{\rho,\alpha}(A) NJ_{\rho,\alpha}(A),
$$
we have
$$
NI_{\rho,\alpha}(A) \leq  \left( \frac{(\lambda_n^{\alpha}-\lambda_1^{\alpha})(\lambda_n^{1-\alpha}-\lambda_1^{1-\alpha})}{(\lambda_n^{\alpha}+\lambda_1^{\alpha})(\lambda_n^{1-\alpha}+\lambda_1^{1-\alpha})} \right)^{1/2} NU_{\rho,\alpha}(A).
$$
Also since
$$
NI_{\rho,\alpha}(A) NJ_{\rho,\alpha}(A) \leq  \frac{(\lambda_n^{\alpha}-\lambda_1^{\alpha})(\lambda_n^{1-\alpha}-\lambda_1^{1-\alpha})}{(\lambda_n^{\alpha}+\lambda_1^{\alpha})(\lambda_n^{1-\alpha}+\lambda_1^{1-\alpha})} NJ_{\rho,\alpha}(A)^2,
$$
we have
$$
NU_{\rho,\alpha}(A) \leq  \left( \frac{(\lambda_n^{\alpha}-\lambda_1^{\alpha})(\lambda_n^{1-\alpha}-\lambda_1^{1-\alpha})}{(\lambda_n^{\alpha}+\lambda_1^{\alpha})(\lambda_n^{1-\alpha}+\lambda_1^{1-\alpha})} \right)^{1/2} NJ_{\rho,\alpha}(A).
$$
Then we have
\begin{eqnarray*}
&   & \left( \frac{(\lambda_n^{\alpha}+\lambda_1^{\alpha})(\lambda_n^{1-\alpha}+\lambda_1^{1-\alpha})}{(\lambda_n^{\alpha}-\lambda_1^{\alpha})(\lambda_n^{1-\alpha}-\lambda_1^{1-\alpha})} \right)^{1/2}  NI_{\rho,\alpha}(A) \\
& \leq & NU_{\rho,\alpha}(A) \leq \left( \frac{(\lambda_n^{\alpha}-\lambda_1^{\alpha})(\lambda_n^{1-\alpha}-\lambda_1^{1-\alpha})}{(\lambda_n^{\alpha}+\lambda_1^{\alpha})(\lambda_n^{1-\alpha}+\lambda_1^{1-\alpha})} \right)^{1/2} NJ_{\rho,\alpha}(A). 
\end{eqnarray*}
(3)~~We show the relation between $NU_{\rho,\alpha}(A)$ and $NV_{\rho}(A)$. We have 
\begin{eqnarray}
&   & NU_{\rho,\alpha}(A)^2 \nonumber \\
& = & \frac{1}{4} \left\{ (Tr[\rho (A_0^*A_0+A_0A_0^*)])^2 -(Tr[\rho^{\alpha}A_0^* \rho^{1-\alpha}A_0 + \rho^{\alpha}A_0 \rho^{1-\alpha}A_0^*] )^2 \right\} \nonumber \\
& = & \frac{1}{4} (Tr[\rho (A_0^*A_0+A_0A_0^*)])^2 \left\{ 1- \left( \frac{Tr[\rho^{\alpha}A_0^* \rho^{1-\alpha}A_0 + \rho^{\alpha}A_0 \rho^{1-\alpha}A_0^*]}{Tr[\rho (A_0^*A_0+A_0A_0^*)]}\right)^2 \right\}. \label{eq:num2.3}
\end{eqnarray}
Since
$$
\min \left\{ \frac{t^{\alpha}+t^{1-\alpha}}{t+1} ; \frac{\lambda_1}{\lambda_n} \leq t \leq \frac{\lambda_n}{\lambda_1} \right\} = \frac{\lambda_n^{\alpha}\lambda_1^{1-\alpha}+\lambda_n^{1-\alpha}\lambda_1^{\alpha}}{\lambda_n+\lambda_1},
$$
we have
\begin{eqnarray*}
&   & Tr[\rho^{\alpha}A_0^* \rho^{1-\alpha}A_0 + \rho^{\alpha}A_0 \rho^{1-\alpha}A_0^*] \\
& = & \sum_{i,j} (\lambda_i^{\alpha}\lambda_j^{1-\alpha}+\lambda_i^{1-\alpha}\lambda_j^{\alpha})|a_{ij}|^2 \\
& = & \sum_{i,j} \frac{\lambda_i^{\alpha}\lambda_j^{1-\alpha}+\lambda_i^{1-\alpha}\lambda_j^{\alpha}}{\lambda_i+\lambda_j} (\lambda_i+\lambda_j)|a_{ij}|^2 \\
& \geq & \min_{i,j} \frac{\lambda_i^{\alpha}\lambda_j^{1-\alpha}+\lambda_i^{1-\alpha}\lambda_j^{\alpha}}{\lambda_i+\lambda_j} \sum_{i,j} (\lambda_i+\lambda_j)|a_{ij}|^2 \\
& = & \frac{\lambda_n^{\alpha}\lambda_1^{1-\alpha}+\lambda_n^{1-\alpha}\lambda_1^{\alpha}}{\lambda_n+\lambda_1} Tr[\rho(A_0^*A_0+A_0A_0^*)].
\end{eqnarray*}
Then since 
$$
\frac{Tr[\rho^{\alpha}A_0^* \rho^{1-\alpha}A_0 + \rho^{\alpha}A_0 \rho^{1-\alpha}A_0^*]}{Tr[\rho (A_0^*A_0+A_0A_0^*)]} \geq \frac{\lambda_n^{\alpha}\lambda_1^{1-\alpha}+\lambda_n^{1-\alpha}\lambda_1^{\alpha}}{\lambda_n+\lambda_1},
$$
we have the following by (\ref{eq:num2.3}).
\begin{eqnarray*}
&   & NU_{\rho,\alpha}(A)^2 \\
& \leq & \frac{(\lambda_n+\lambda_1)^2-(\lambda_n^{\alpha}\lambda_1^{1-\alpha}+\lambda_n^{1-\alpha}\lambda_1^{\alpha})^2}{(\lambda_m+\lambda_1)^2} NV_{\rho}(A)^2 \\
& = & \frac{ (\lambda_n^{\alpha}-\lambda_1^{\alpha})(\lambda_n^{\alpha}+\lambda_1^{\alpha})(\lambda_n^{1-\alpha}-\lambda_1^{1-\alpha})(\lambda_n^{1-\alpha}+\lambda_1^{1-\alpha})}{(\lambda_n+\lambda_1)^2} NV_{\rho}(A)^2.
\end{eqnarray*}
Therefore we get the result.
\ \hfill $\Box$

\vspace{0.5cm}
\begin{remark}
We rewrite the results of Theorem \ref{th:theorem2.1}.
\begin{eqnarray*}
NI_{\rho,\alpha}(A) & \leq & \left( \frac{(\lambda_n^{\alpha}-\lambda_1^{\alpha})(\lambda_n^{1-\alpha}-\lambda_1^{1-\alpha})}{(\lambda_n^{\alpha}+\lambda_1^{\alpha})(\lambda_n^{1-\alpha}+\lambda_1^{1-\alpha})} \right)^{1/2} NU_{\rho,\alpha}(A) \\
& \leq & \frac{(\lambda_n^{\alpha}-\lambda_1^{\alpha})(\lambda_n^{1-\alpha}-\lambda_1^{1-\alpha})}{\lambda_n+\lambda_1} NV_{\rho}(A) \\
& \leq & \frac{(\lambda_n^{\alpha}-\lambda_1^{\alpha})(\lambda_n^{1-\alpha}-\lambda_1^{1-\alpha})}{(\lambda_n^{\alpha}+\lambda_1^{\alpha})(\lambda_n^{1-\alpha}+\lambda_1^{1-\alpha})} NJ_{\rho,\alpha}(A).
\end{eqnarray*}
\label{re:remark2.1}
\end{remark}


\section{Refined Non-Hermite Wigner-Yanase-Dyson Skew Information}

\begin{theorem}
Let $\rho$ be a density matrix with eigenvalues $0 < \lambda_1 \leq \lambda_2 \leq \cdots \leq \lambda_n$ and $\lambda_1 \neq \lambda_n$ and $f(x)$ be non-constant continuous function. Then 
\begin{description}
\item[(1)]  $\displaystyle{\frac{\lambda_n+\lambda_1}{(\lambda_n^{\alpha}-\lambda_1^{\alpha})(\lambda_n^{1-\alpha}-\lambda_1^{1-\alpha})} NI_{\rho,\alpha}(A) + \frac{|Tr[\rho f(\rho)_0 A_0]|^2}{Tr[\rho f(\rho)_0^2]} \leq NV_{\rho}(A)}$. 
\item[(2)]  $\displaystyle{\frac{(\lambda_n^{\alpha}+\lambda_1^{\alpha})(\lambda_n^{1-\alpha}+\lambda_1^{1-\alpha})}{(\lambda_n^{\alpha}-\lambda_1^{\alpha})(\lambda_n^{1-\alpha}-\lambda_1^{1-\alpha})} NI_{\rho,\alpha}(A) + 2\frac{|Tr[\rho f(\rho)_0 A_0]|^2}{Tr[\rho f(\rho)_0^2]} \leq NJ_{\rho,\alpha}(A) }$.
\item[(3)]  $\displaystyle{ \left\{ \frac{(\lambda_n^{\alpha}+\lambda_1^{\alpha})(\lambda_n^{1-\alpha}+\lambda_1^{1-\alpha})}{(\lambda_n^{\alpha}-\lambda_1^{\alpha})(\lambda_n^{1-\alpha}-\lambda_1^{1-\alpha})} NI_{\rho,\alpha}(A)^2 +2\frac{|Tr[\rho f(\rho)_0 A_0]|^2}{Tr[\rho f(\rho)_0^2]}NI_{\rho,\alpha}(A) \right\}^{1/2}  \leq NU_{\rho,\alpha}(A)}$.
\end{description}
\label{th:theorem3.1}
\end{theorem}

\noindent
Proof of Theorem \ref{th:theorem3.1}~~We remark that 
$$
A_0 = A-Tr[\rho A]I,~~A_0^* = A^*-\overline{Tr[\rho A]}I, 
$$
$$
f(\rho)_0 = f(\rho)-Tr[\rho f(\rho)]I,~~f(\rho)_0^* = f(\rho)-Tr[\rho f(\rho)]I.
$$
$$
(A+tf(\rho))_0 = A_0+t f(\rho)_0,~~~(A+t f(\rho))_0^* = (A_0+t f(\rho)_0)^* = A_0^*+\bar{t} f(\rho)_0.
$$
Since $[\rho^{1-\alpha},A_0] = [\rho^{1-\alpha},A_0+tf(\rho)_0]$ and $[\rho^{\alpha},A_0^*] = [\rho^{\alpha},A_0^*+\bar{t}f(\rho)_0)]$ for $t \in \mathbb{C}$, we have $NI_{\rho,\alpha}(A) = NI_{\rho,\alpha}(A+tf(\rho))$.

\vspace{0.3cm}
\noindent
(1)~~By Theorem \ref{th:theorem2.1} (1) 
$$
\frac{\lambda_n+\lambda_1}{(\lambda_n^{\alpha}-\lambda_1^{\alpha})(\lambda_n^{1-\alpha}-\lambda_1^{1-\alpha})} NI_{\rho,\alpha}(A+tf(\rho)) \leq NV_{\rho}(A+tf(\rho)).
$$
\begin{eqnarray*}
&   & NV_{\rho}(A+tf(\rho)) \\
& = & \frac{1}{2} Tr[\rho \{ (A_0^*+\bar{t}f(\rho)_0)(A_0+tf(\rho)_0)+(A_0+tf(\rho)_0)(A_0^*+\bar{t}f(\rho)_0) \}] \\
& = & \frac{1}{2} Tr[\rho \{ A_0^*A_0+tA_0^*f(\rho)_0+\bar{t}f(\rho)_0A_0+|t|^2f(\rho)_0^2+A_0A_0^*+tf(\rho)_0A_0^*+\bar{t}A_0f(\rho)_0+|t|^2f(\rho)_0^2 \}] \\
& = & \frac{1}{2} Tr[\rho(A_0^*A_0+A_0A_0^*)]+|t|^2Tr[\rho f(\rho)_0^2]+t Tr[\rho f(\rho)_0 A_0^*]+\bar{t} Tr[\rho f(\rho)_0 A_0] \\
& = & \frac{1}{2} Tr[\rho(A_0^*A_0+A_0A_0^*)]+|t|^2Tr[\rho f(\rho)_0^2] + 2 Re \{ tTr[\rho f(\rho)_0A_0^*] \}.
\end{eqnarray*}
By putting $t = -\frac{\overline{Tr[\rho f(\rho)_0A_0^*]}}{Tr[\rho f(\rho)_0^2]}$, we have 
$$
NV_{\rho}(A+tf(\rho)) = NV_{\rho}(A) -\frac{|Tr[\rho f(\rho)_0A_0^*]|^2}{Tr[\rho f(\rho)_0^2]}.
$$
Hence 
$$
\frac{\lambda_n+\lambda_1}{(\lambda_n^{\alpha}-\lambda_1^{\alpha})(\lambda_n^{1-\alpha}-\lambda_1^{1-\alpha})} NI_{\rho,\alpha}(A) \leq NV_{\rho}(A)-\frac{|Tr[\rho f(\rho)_0A_0^*]|^2}{Tr[\rho f(\rho)_0^2]}. 
$$
Since $|Tr[\rho f(\rho)_0A_0^*]| = |Tr[\rho f(\rho)_0A_0]|$, we have 
$$
\frac{\lambda_n+\lambda_1}{(\lambda_n^{\alpha}-\lambda_1^{\alpha})(\lambda_n^{1-\alpha}-\lambda_1^{1-\alpha})} NI_{\rho,\alpha}(A) +\frac{|Tr[\rho f(\rho)_0A_0]|^2}{Tr[\rho f(\rho)_0^2]} \leq NV_{\rho}(A).
$$

\vspace{0.5cm}
\noindent
(2)~~By Theorem \ref{th:theorem2.1} (1) 
$$
\frac{(\lambda_n^{\alpha}+\lambda_1^{\alpha})(\lambda_n^{1+\alpha}+\lambda_1^{1-\alpha})}{(\lambda_n^{\alpha}-\lambda_1^{\alpha})(\lambda_n^{1-\alpha}-\lambda_1^{1-\alpha})} NI_{\rho,\alpha}(A+tf(\rho)) \leq NJ_{\rho,\alpha}(A+tf(\rho)).
$$
\begin{eqnarray*}
&   & NJ_{\rho,\alpha}(A+tf(\rho)) \\
& = & NV_{\rho}(A+ tf(\rho)) + \frac{1}{2} Tr[\rho^{\alpha}(A_0^*+\bar{t}f(\rho)_0) \rho^{1-\alpha} (A_0+tf(\rho)_0)] \\
&   & +\frac{1}{2} Tr[\rho^{\alpha}(A_0+tf(\rho)_0) \rho^{1-\alpha} (A_0^*+\bar{t}f(\rho)_0)] \\
& = & \frac{1}{2} Tr[\rho(A_0^*A_0+A_0A_0^*)]+|t|^2Tr[\rho f(\rho)_0^2] + 2 Re \{ tTr[\rho f(\rho)_0A_0^*] \} \\
&   & + \frac{1}{2} Tr[\rho^{\alpha}(A_0^*+\bar{t}f(\rho)_0) \rho^{1-\alpha} (A_0+tf(\rho)_0)] +\frac{1}{2} Tr[\rho^{\alpha}(A_0+tf(\rho)_0) \rho^{1-\alpha} (A_0^*+\bar{t}f(\rho)_0)] \\
& = & \frac{1}{2} Tr[\rho(A_0^*A_0+A_0A_0^*)]+\frac{1}{2} Tr[\rho^{\alpha}A_0^*\rho^{1-\alpha}A_0+\rho^{\alpha}A_0\rho^{1-\alpha}A_0^*] \\
&   & +2|t|^2 Tr[\rho f(\rho)_0^2] + 4 Re \{ tTr[\rho f(\rho)_0A_0^*] \}.
\end{eqnarray*}
By putting $t = -\frac{\overline{Tr[\rho f(\rho)_0 A_0^*]}}{Tr[\rho f(\rho)_0^2]}$, we have 
$$
NJ_{\rho,\alpha}(A+tI) = NJ_{\rho,\alpha}(A) -2\frac{|Tr[\rho f(\rho)_0A_0^*]|^2}{Tr[\rho f(\rho)_0^2]}. 
$$
Hence 
$$
\frac{(\lambda_n^{\alpha}+\lambda_1^{\alpha})(\lambda_n^{1+\alpha}+\lambda_1^{1-\alpha})}{(\lambda_n^{\alpha}-\lambda_1^{\alpha})(\lambda_n^{1-\alpha}-\lambda_1^{1-\alpha})} NI_{\rho,\alpha}(A)+ 2\frac{|Tr[\rho f(\rho)_0A_0^*]|^2}{Tr[\rho f(\rho)_0^2]} \leq NJ_{\rho,\alpha}(A).
$$
Since $|Tr[\rho f(\rho)_0A_0^*]| = |Tr[\rho f(\rho)_0A_0]|$, we have 
$$
\frac{(\lambda_n^{\alpha}+\lambda_1^{\alpha})(\lambda_n^{1+\alpha}+\lambda_1^{1-\alpha})}{(\lambda_n^{\alpha}-\lambda_1^{\alpha})(\lambda_n^{1-\alpha}-\lambda_1^{1-\alpha})}NI_{\rho,\alpha}(A)+ 2\frac{|Tr[\rho f(\rho)_0A_0]|^2}{Tr[\rho f(\rho)_0^2]} \leq NJ_{\rho,\alpha}(A).
$$

\vspace{0.5cm}
\noindent
(3)~~By (2) we have 
$$
\frac{(\lambda_n^{\alpha}+\lambda_1^{\alpha})(\lambda_n^{1-\alpha}+\lambda_1^{1-\alpha})}{(\lambda_n^{\alpha}-\lambda_1^{\alpha})(\lambda_n^{1-\alpha}-\lambda_1^{1-\alpha})} NI_{\rho,\alpha}(A)^2 +2\frac{|Tr[\rho f(\rho)_0A_0]|^2}{Tr[\rho f(\rho)_0^2]} NI_{\rho,\alpha}(A) \leq NU_{\rho,\alpha}(A)^2. 
$$
Then we obtain the result. \ \hfill $\Box$

\vspace{0.5cm}
\begin{remark}
In the two dimensional case, the equalities in (1), (2) and (3) in Theorem \ref{th:theorem3.1} hold.  Let $\rho = \left( \begin{array}{cc}
\lambda_1 & 0 \\
0 & \lambda_2
\end{array} \right)$ and $A = \left( \begin{array}{cc}
a & c \\
d & b
\end{array} \right)$. Then 
$$
NI_{\rho,\alpha}(A) =  \frac{1}{2}(\lambda_1^{\alpha}-\lambda_2^{\alpha})(\lambda_1^{1-\alpha}-\lambda_2^{1-\alpha})(|c|^2+|d|^2),
$$
$$
NJ_{\rho,\alpha}(A) = 2\lambda_1\lambda_2|a-b|^2+\frac{1}{2}(\lambda_1^{\alpha}+\lambda_2^{\alpha})(\lambda_1^{1-\alpha}+\lambda_2^{1-\alpha})(|c|^2+|d|^2),
$$
$$
NV_{\rho}(A) = \lambda_1\lambda_2|a-b|^2+\frac{1}{2}(|c|^2+|d|^2),
$$
Since $Tr[\rho f(\rho)_0^2] = \lambda_1\lambda_2(f(\lambda_1)-f(\lambda_2))^2$ and $Tr[\rho f(\rho)_0A_0] = \lambda_1\lambda_2(a-b)(f(\lambda_1)-f(\lambda_2))$, we have 
\begin{eqnarray*}
&   & \frac{\lambda_1+\lambda_2}{(\lambda_1^{\alpha}-\lambda_2^{\alpha})(\lambda_1^{1-\alpha}-\lambda_2^{1-\alpha})}NI_{\rho,\alpha}(A) + \frac{|Tr[\rho f(\rho)_0 A_0]|^2}{Tr[\rho f(\rho)_0^2]} \\
& = & \frac{1}{(\lambda_1^{\alpha}-\lambda_2^{\alpha})(\lambda_1^{1-\alpha}-\lambda_2^{1-\alpha})} \frac{1}{2}(\lambda_1^{\alpha}-\lambda_2^{\alpha})(\lambda_1^{1-\alpha}-\lambda_2^{1-\alpha})(|c|^2+|d|^2) + \lambda_1\lambda_2|a-b|^2 \\
& = & \frac{1}{2}(|c|^2+|d|^2) + \lambda_1\lambda_2|a-b|^2 = NV_{\rho}(A).
\end{eqnarray*}
In (2) and (3) we obtain the similar results.
\label{re:remark3.1}
\end{remark} 


\end{document}